\documentclass[12pt]{article}
\usepackage{amssymb,amsfonts,amsmath, psfrag,eepic,colordvi, graphics}
\topskip  
\numberwithin{equation}{section}
\newtheorem{theorem}{Theorem}[section]

\newtheorem{conjecture}[theorem]{Conjecture}

\newtheorem{lemma}[theorem]{Lemma}

\DeclareMathOperator{\as}{asc}

\begin{document}
\parskip 6pt

\pagenumbering{arabic}
\def\sof{\hfill\rule{2mm}{2mm}}
\def\ls{\leq}
\def\gs{\geq}
\def\SS{\mathcal S}
\def\qq{{\bold q}}
\def\MM{\mathcal M}
\def\TT{\mathcal T}
\def\EE{\mathcal E}
\def\lsp{\mbox{lsp}}
\def\rsp{\mbox{rsp}}
\def\pf{\noindent {\it Proof.} }
\def\mp{\mbox{pyramid}}
\def\mb{\mbox{block}}
\def\mc{\mbox{cross}}
\def\qed{\hfill \rule{4pt}{7pt}}
\def\block{\hfill \rule{5pt}{5pt}}

\begin{center}
{\Large\bf Bijections for inversion sequences, ascent sequences and  3-nonnesting set partitions }
\vskip 6mm
\end{center}

\begin{center}
{\small
Sherry H. F. Yan \\[2mm]
 Department of Mathematics, Zhejiang Normal University, Jinhua
321004, P.R. China
\\[2mm]
  huifangyan@hotmail.com
 \\[0pt]
}
\end{center}

\noindent {\bf Abstract.} Set partitions avoiding $k$-crossing and $k$-nesting have been extensively studied   from the aspects of both combinatorics and mathematical biology.     By using
the generating tree technique,   the obstinate kernel method and Zeilberger's algorithm, Lin  confirmed    a  conjecture due independently to the author and Martinez-Savage that
asserts inversion sequences    with no weakly decreasing subsequence of length 3 and enhanced
3-nonnesting partitions  have the same cardinality.   In this paper, we provide a bijective proof of this conjecture. Our bijection also enables us to provide a new bijective  proof of a conjecture posed by Duncan and Steingr\'{\i}msson,  which was  proved by the author via  an intermediate structure of growth diagrams for $01$-fillings of Ferrers shapes.




\section{Introduction}

Set partitions avoiding $k$-crossing and $k$-nesting have been extensively studied  from the aspects of both combinatorics and mathematical biology; see $\cite{chen,chen3,kra}$ and the references therein.
The objective of this paper is to  provide a bijective proof of a conjecture due independently to the author \cite{Yan} and Martinez-Savage \cite{Mar}, which was recently confirmed by Lin \cite{Lin}    using  the generating tree technique, the obstinate kernel method \cite{melon1}  and Zeilberger's algorithm \cite{Zei}. Our bijection also enables us to provide a new bijective  proof of a conjecture posed by Duncan and Steingr\'{\i}msson \cite{Duncan},  which   was   proved by the author \cite{Yan} via  an intermediate structure of growth diagrams for $01$-fillings of Ferrers shapes  \cite{kra} and \cite{van}.
Let us first give an overview of the notation and terminology.

A sequence $x=x_1 x_2  \cdots x_n $ is said to be an  {\em inversion
sequence} of length $n$  if it satisfies   $0\leq x_i<
i$ for all $1\leq i\leq n$.  Inversion sequences of length $n$ are in easy bijection with permutations of length $n$. An inversion sequence $x_1x_2\ldots x_n$ can be obtained from any permutation $\pi=\pi_1\pi_2\ldots \pi_n$ by setting $x_i=\{j\mid j<i \,\,\mbox{ and }\,\, \pi_{j}>\pi_{i}\}$.

Given a sequence of integers   $x=x_1x_2\cdots x_n$, we say
that the sequence $x$ has an {\em ascent} at position $i$ if
$x_i<x_{i+1}$. The number of ascents of $x$ is denoted by $\as(x)$.
A sequence $x=x_1x_2\cdots x_n$ is said to be an {\em ascent
sequence of length $n$} if it satisfies $x_1=0$ and $0\leq x_i\leq
\as(x_1x_2\cdots x_{i-1})+1$ for all $2\leq i\leq n$. Ascent
sequences were introduced by Bousquet-M$\acute{e}$lou et al.
\cite{melon} in their study of $(2+2)$-free posets, which  are  closely connected to unlabeled  $(2+2)$-free posets,  permutations avoiding a certain pattern,
 and a class of involutions
introduced by Stoimenow \cite{Stoi}.
  We call an ascent sequence with no two consecutive equal
entries a {\em primitive} ascent sequence.

Pattern avoiding permutations have been extensively studied  over last decade.
For a thorough summary of the
current status of research,
see B\'{o}na's  book \cite{bona} and Kitaev's book \cite{kitaev}. Analogous to pattern avoidance on permutations,  Corteel-Martinez-Savage-Weselcouch \cite{Cor}  and Mansour-Shattuck \cite{man2} initiated the study of inversion  sequences avoiding certain pattern. Pattern avoiding inversion sequences   are closely  related to Catalan numbers,  large Schr$\ddot{o}$der numbers, Euler numbers and Baxter  numbers (see \cite{Cor}, \cite{Kim}, \cite{Lin} and \cite{man2}).  In their paper  \cite{Duncan}, Duncan and Steingr\'{\i}msson
studied  ascent sequences avoiding certain patterns.      Further  results on the
 enumeration of pattern-avoiding ascent sequences could be found in  \cite{chen1, man,Yan}.
   By using
the generating tree technique,  the  obstinate kernel method and Zeilberger's algorithm, Lin \cite{Lin}  confirmed  the following     conjecture proposed by Martinez-Savage \cite{Mar}.

\begin{conjecture}{ \upshape    (  Martinez-Savage  \cite{Mar} )}\label{con}
    Inversion sequences of length $n$  and with  no  weakly decreasing subsequence of length $3$ are equinumerous with enhanced  3-nonnesting (3-noncrossing) set partitions of
$[n]$.
\end{conjecture}
As remarked by Lin \cite{Lin}, this conjecture has already been proposed by Yan \cite{Yan} in the the course of confirming  the
following  conjecture  posed by Duncan and Steingr\'{\i}msson \cite{Duncan}.

 \begin{conjecture}{ \upshape    ( See \cite{Duncan}, Conjecture 3.3)}\label{Yan}
 Ascent sequences of length $n$ and with  no   decreasing subsequence of length $3$ are equinumerous with  3-nonnesting (3-noncrossing) set partitions of
$[n]$.
 \end{conjecture}

 Recall that a subsequence $x_{i_1}x_{i_2}\ldots x_{i_k} $ of a  sequence $x=x_1x_2\ldots x_n$ is said to be  {\em decreasing }  if $i_1<i_2<\ldots  <i_k$ and $x_{i_1}> x_{i_2}>\ldots> x_{i_k}$ and to be  {\em weakly decreasing}
if $i_1<i_2<\ldots  <i_k$ and $x_{i_1}\geq  x_{i_2}\geq \ldots\geq  x_{i_k}$.
Denote by $\mathcal{A}_k(n)$ and $\mathcal{PA}_k(n)$ the set of ordinary and primitive    ascent sequences of length $n$ and  with no decreasing subsequence  of length $k$, respectively. Let $\mathcal{I}_k(n)$ denote  the set of inversion sequences of length $n$ and with no weakly decreasing sequences of length $k$.

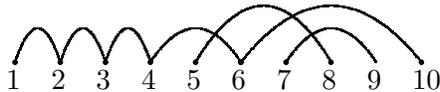
\begin{figure}[h!]
\begin{center}
\begin{picture}(50,15)
\setlength{\unitlength}{6mm} \linethickness{0.4pt}
 \put(-0.2,0.4){\small$1$}\put(0.8,0.4){\small$2$}
\put(1.8,0.4){\small$3$}\put(2.8,0.4){\small$4$}
 \put(3.8,0.4){\small$5$} \put(4.8,0.4){\small$6$}
\put(5.8,0.4){\small$7$}\put(6.8,0.4){\small$8$}\put(7.8,0.4){\small$9$}\put(8.8,0.4){\small$10$}
\qbezier(0,1)(0.5,2.5)(1,1)\qbezier(1,1)(1.5,2.5)(2,1)\qbezier(2,1)(2.5,2.5)(3,1)\qbezier(3,1)(4,2.5)(5,1)
\qbezier(4,1)(5.5,3.5)(7,1)\qbezier(5,1)(7,3.5)(9,1)\qbezier(6,1)(7,2.5)(8,1)
\put(5,1){\circle*{0.1}} \put(0,1){\circle*{0.1}}\put(1,1){\circle*{0.1}}\put(2,1){\circle*{0.1}}
\put(3,1){\circle*{0.1}}\put(4,1){\circle*{0.1}}\put(5,1){\circle*{0.1}}\put(6,1){\circle*{0.1}}\put(7,1){\circle*{0.1}}
\put(9,1){\circle*{0.1}}

\end{picture}
\vspace{-15pt}
\end{center}
\caption{ The linear representation of  a set partition $\pi=\{\{1,2,3,4, 6,10\}, \{5,8\}, \{7,9\}\}$ .} \label{linear}
\end{figure}

  A set partition $P$ of $[n] = \{1,2,\cdots, n\}$ can be represented
by a diagram with vertices drawn on a horizontal line in increasing order. For a block
$B$ of $P$, we write the elements of $B$ in increasing order. Suppose that $B=\{i_1, i_2, \cdots, i_k\}$.
Then we draw an arc from $i_1$ to $i_2$, an arc from $i_2$ to $i_3$, and so on. Such a diagram is called
the {\em linear representation} of $P$, see Figure \ref{linear} for  example. The {\em enhanced }  representation of $P$ is defined to the union of the standard representation of $P$ and the set of loops $(i,i)$, where $i$ ranges over all the singleton blocks $\{i\}$ of $P$. Then one defines a {\em $k$-crossing} of a set partition to be a subset $\{(i_1, j_1), (i_2, j_2), \cdots, (i_k, j_k)\}$ of its  linear representation where   $i_1<i_2<\cdots <i_k< j_1<j_2<\cdots<j_k$, and   an {\em enhanced $k$-crossing } of a set partition to be a subset $\{(i_1, j_1), (i_2, j_2), \cdots, (i_k, j_k)\}$ of its enhanced representation where   $i_1<i_2<\cdots <i_k\leq j_1<j_2<\cdots<j_k$.  A partition without any  (enhanced) $k$-crossings  is said to be {\em  (enhanced )$k$-noncrossing}.
Similarily,   a  {\em $k$-crossing } is  defined to be a subset $\{(i_1, j_1), (i_2, j_2), \cdots, (i_k, j_k)\}$ of its enhanced representation where   $i_1<i_2<\cdots <i_k<j_k<j_{k-1}<\cdots<j_1$, and  an {\em enhanced  $k$-nesting} is  defined to be a subset $\{(i_1, j_1), (i_2, j_2), \cdots, (i_k, j_k)\}$ of its enhanced representation where   $i_1<i_2<\cdots <i_k\leq j_k<j_{k-1}<\cdots<j_1$.
 A set partition without any (enhanced )   $k$-nestings  is said to be {\em  (enhanced ) $k$-nonnesting}. Chen et al. \cite{chen} proved that   (enhanced) $k$-nonnesting set partitions of $[n]$ are equinumerous with  (enhanced) $k$-noncrossing set
partitions of $[n]$ bijectively using hesitating   tableaux as an intermediate object. Denote by $\mathcal{C}_k(n)$ and $\mathcal{E}_k(n)$ the set of  ordinary and enhanced  $k$-nonnesting   set partitions of $[n]$, respectively.

\section{ Bijective proof of Conjecture \ref{con} }
In this section, we shall  provide a bijective proof of Conjecture \ref{con} by showing that inversion sequences  of length $n$ and    with no weakly decreasing subsequence of length 3  are in bijection with  enhanced
3-nonnesting partitions of $[n]$.
 To this end,  we recall some necessary   notation and terminology.

A {\em triangular shape} of order $n$ is  the
left-justified array of  ${n+1\choose 2}$ squares    in which the $i$th row contains exactly $i$ squares.  Let $\Delta_n$ be the triangular shape of order $n$. In a triangular shape,  we number rows from top to bottom and columns from left to right and identify squares  using matrix coordinate. The $i$th row (column) is called row (column) $i$.  For example,  the square  in
the first row and second column is numbered $(1, 2)$.

\begin{figure}[h!]
\begin{center}
\begin{picture}(50,20)
\setlength{\unitlength}{0.6mm}
\put(0,0){\framebox(5,5) } \put(5,0){\framebox(5,5)    }  \put(10,0){\framebox(5,5) {$\bullet$}  }
 \put(15,0){\framebox(5,5) } \put(20,0){\framebox(5,5) }  \put(25,0){\framebox(5,5) }

\put(0,5){\framebox(5,5) } \put(5,5){\framebox(5,5)    }  \put(10,5){\framebox(5,5)  }
 \put(15,5){\framebox(5,5) {$\bullet$}  }
\put(20,5){\framebox(5,5) }

\put(0,10){\framebox(5,5) } \put(5,10){\framebox(5,5)  {$\bullet$} }  \put(10,10){\framebox(5,5)  }
 \put(15,10){\framebox(5,5)  {$\bullet$}}

 \put(0,15){\framebox(5,5) } \put(5,15){\framebox(5,5)  }  \put(10,15){\framebox(5,5)  }

 \put(0,20){\framebox(5,5)  } \put(5,20){\framebox(5,5)  {$\bullet$} }

  \put(0,25){\framebox(5,5)  {$\bullet$} }

  \end{picture}
\end{center}
\caption{An  example of a  $01$-filling of a triangular shape of order $6$.  }\label{afilling}
\end{figure}
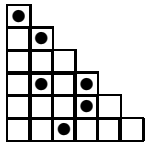

  A  {\em $01$-filling} of a triangular shape   $\Delta_n$ is obtained by filling the squares  of $\Delta_n$ with $1's$ and $0's$, see Figure \ref{afilling} for    example,  where we represent a $1$ by a $\bullet$ and suppress the $0$'s. A  $01$-filling of a triangular shape is said to be {\em valid} if every row contains at most one $1$. A row (column)  of a $01$-filling  is said to a {\em zero}  if   all the squares at this row (  column) are filled with $0's$.
A {\em NE-chain} of a $01$-filling is a sequence of $1's$ such that any $1$ is strictly above and weakly to the right of the preceding $1$ in the sequence.  For example, in Figure \ref{afilling},  the  sequence of $1's$ lying in the squares  $(6,3)$, $(5,4,)$ and  $(4,4)$    form a NE-chain of length $3$.

 An inversion sequence $x_1x_2\ldots x_n$   can  be encoded by a $01$-filling of $\Delta_n$ in which the square $(i,x_i+1)$ is filled with a $1$ for all $1\leq i\leq n$ and all the other squares are filled with $0's$.  It is easily seen that a weakly decreasing sequence of length $k$ corresponds to a NE-chain of length $k$. Denote by $\mathcal{M}_k(n)$ the set of $01$-fillings of $\Delta_{n}$ with the property  that
     every row   contains exactly one $1$ and
    there is no NE-chain of length $k$.

\begin{theorem}\label{th1}
There is a one-to-one correspondence  between the set $\mathcal{I}_k(n)$ and the set $\mathcal{M}_k(n)$.
\end{theorem}

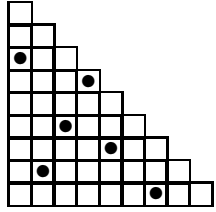
\begin{figure}[h!]
\begin{center}
\begin{picture}(200,30)
\setlength{\unitlength}{0.6mm}
\put(60,0){\framebox(5,5) } \put(65,0){\framebox(5,5)  }  \put(70,0){\framebox(5,5)  }
 \put(75,0){\framebox(5,5) }  \put(80,0){\framebox(5,5) }  \put(85,0){\framebox(5,5) }
  \put(90,0){\framebox(5,5){$\bullet$}  }  \put(95,0){\framebox(5,5)  }  \put(100,0){\framebox(5,5)  }

 \put(60,5){\framebox(5,5)  } \put(65,5){\framebox(5,5)  {$\bullet$}}  \put(70,5){\framebox(5,5)  }
 \put(75,5){\framebox(5,5) }  \put(80,5){\framebox(5,5) }  \put(85,5){\framebox(5,5)  }
  \put(90,5){\framebox(5,5) }  \put(95,5){\framebox(5,5)  }

 \put(60,10){\framebox(5,5)  } \put(65,10){\framebox(5,5)  }  \put(70,10){\framebox(5,5)  }
 \put(75,10){\framebox(5,5) }  \put(80,10){\framebox(5,5) {$\bullet$} }  \put(85,10){\framebox(5,5)  }
  \put(90,10){\framebox(5,5) }

  \put(60,15){\framebox(5,5)  } \put(65,15){\framebox(5,5)  }  \put(70,15){\framebox(5,5){$\bullet$}  }
 \put(75,15){\framebox(5,5)}  \put(80,15){\framebox(5,5) }  \put(85,15){\framebox(5,5)  }

\put(60,20){\framebox(5,5)  } \put(65,20){\framebox(5,5)  }  \put(70,20){\framebox(5,5)  }
 \put(75,20){\framebox(5,5)}  \put(80,20){\framebox(5,5) }

\put(60,25){\framebox(5,5) } \put(65,25){\framebox(5,5)  }  \put(70,25){\framebox(5,5)  }
 \put(75,25){\framebox(5,5){$\bullet$} }

\put(60,30){\framebox(5,5){$\bullet$}  } \put(65,30){\framebox(5,5)   }  \put(70,30){\framebox(5,5)  }

\put(60,35){\framebox(5,5)  } \put(65,35){\framebox(5,5)  }

\put(60,40){\framebox(5,5)  }

  \end{picture}
\end{center}
\caption{ A set partition $\pi=\{\{1, 3,6 \},  \{2, 8\}, \{4\}, \{5,7,9\}\}$   and its corresponding   $01$-filling.}\label{filling}
\end{figure}

In his paper \cite{kra},     Krattenthaler  established      a bijection  between set partitions of $[n]$ and $01$-fillings of $\Delta_{n}$ in which
   every row and every column contain  at most one $1$, and
   either column $i$ or row $i$ contains  at least one $1$  for all $1\leq i\leq n$.
For the sake of completeness, we give a brief description of this bijection.   Given a set partition $\pi$ of $[n]$, we can get a $01$-filling of $\Delta_{n}$ by putting a $1$ in the square $(j,i)$  if  $(i,j)$ is an arc in its enhanced representation, and, in addition, by putting a   $1$ in the the square $(i,i)$   if $(i,i)$ is a loop in its enhanced representation. The $01$-filling corresponding to the set partition $\pi=\{\{1, 3,6 \}, \{2, 8\}, \{4\}, \{5,7,9\}\}$ is indicated  in Figure \ref{filling}.   From the construction of Krattenthaler's bijection,    an enhanced  $k$-nesting of a set partition  corresponds to a NE-chain of length $k$ in its corresponding $01$-filling.

 Denote by $\mathcal{N}_k(n)$ the set of $01$-fillings of $\Delta_{n}$ satisfying the  following properties:
\begin{itemize}
\item[(a1)]  every row and every column contain  at most one $1$;
\item[(b1)] either column $i$ or row $i$ contains at least one $1$  for all $1\leq i\leq n$;
\item[(c1)] there is no NE-chain of length $k$.
\end{itemize}
From Krattenthaler's bijection, we immediately  get the following result.

\begin{theorem}\label{th2}
   The $01$-fillings of  the set $\mathcal{E}_k(n)$ are in  bijection with the $01$-fillings of  the set $\mathcal{N}_k(n)$.
\end{theorem}

In view of Theorems \ref{th1} and \ref{th2}, in order to prove Conjecture \ref{con}, it suffices to establish a bijection between the set $\mathcal{M}_k(n)$ and the set $\mathcal{N}_k(n)$.  To this end, we define two
   transformations, which will play an essential role in the construction of the bijection.

{\noindent \bf The transformation $\alpha$}\,\, Let $F$ be a valid $01$-filling of $\Delta_n$ without any NE-chain of length $3$. If every column of $F$ contains at most one $1$, we simply define $\alpha(F)=F$. Otherwise, find the leftmost column $i$ which  contains  at least two $1's$. Suppose  the square  $(i,  j)$ is  filled with a  $1$ for some $1\leq j\leq i$.
Assume that  the $1's$ below and  weakly to the left of  the square $(i,i)$ are positioned at the squares $(r_1, c_1)$, $(r_2, c_2), \ldots, (r_m, c_m)$ with $r_1<r_2<\ldots <r_m$. Let $r_0=i, c_0=j$.  Suppose that the topmost $1$ in column $i$ is at row $r_s$.

 If row $r_s+1$ contains a $1$ which is to the right of the square $(r_s, c_s)$, then define $\alpha (F)$ to be the $01$-filling of $\Delta_n$ obtained   from $F$ by the following procedure:
\begin{itemize}
 \item For all $0\leq \ell\leq s$, replace the $1$ at     the square  $(r_\ell, c_\ell)$ with a $0$;
   \item  For all $0\leq \ell<s $,  fill  the square  $(r_{\ell+1}, c_\ell)$  with a $1 $;

   \item Leave all the other squares fixed.
\end{itemize}

Otherwise,  define $\alpha (F)$ to be the $01$-filling of $\Delta_n$ obtained   from $F$ by the following procedure:
\begin{itemize}
 \item For all $0\leq \ell\leq m$, replace the $1$ at    the square   $(r_\ell, c_\ell)$ with  a $0 $;
   \item  For all $0\leq \ell<m $,  fill  the square  $(r_{\ell+1}, c_\ell)$  with a $1$;

   \item Leave all the other squares fixed.
\end{itemize}

  Now we  proceed to   show that the transformation  $\alpha$   has the following  desired properties.
\begin{lemma}\label{alpha0}
In $\alpha(F)$, each column to the left of column  $i$ contains  at most one $1$ and  column $i$  contains exactly one $1$.
\end{lemma}
\pf It is obvious from the selection of column $i$ and  the construction of the transformation $\alpha$. \qed

\begin{lemma}\label{alpha1}
   The filling $\alpha(F)$  is a valid $01$-filling of $\Delta_n$ containing no  NE-chain of length $3$.
 \end{lemma}
\pf According to the construction of the transformation $\alpha$, it is easily seen that $\alpha(F)$  is a valid $01$-filling of $\Delta_n$. Now we proceed to show that $\alpha(F)$ contains no  NE-chain of length $3$. If not, suppose that the $1's$ positioned at the squares $(a_1, b_1)$, $(a_2, b_2)$ and $(a_3, b_3)$ form a NE-chain of length $3$, where $a_1<a_2<a_3$.
Since $F$ has no NE-chain of length $3$, the square $(a_1, b_1)$ must be positioned below and to the right of the square $(r_s,c_s)$.   Suppose that there exists a $1$ at row $r_s+1$ which is to the right of the square $(r_s, c_s)$. From the construction of $\alpha(F)$, it is easy to check that  all the squares below row $r_s$ remain the same as those of $F$. This implies that there is no NE-chain of length $3$ below row $r_s$ in $\alpha(F)$. This contradicts the fact that $(a_1, b_1)$ is below row $r_s$.  Hence, row $r_s+1$ does not contain a $1$ which is to the right of the square $( r_s, c_s)$.  This implies that the square $(a_1,b_1)$ is below and to the right of the square  $(r_{s+1}, c_s)$. Since $F$ has no NE-chain of length $3$, there is no NE-chain of length $2$ below and to the left of the square $(r_{s+1}, c_s)$ in $\alpha(F)$. This yields that  both the square $(a_1, b_1 )$ and the square $(a_2, b_2)$ are positioned     to the right of the square $(r_s, c_s)$.  From the fact that $F$ contains no NE-chain of length $3$, we have $c_s=c_m=i$.  Then the $1's$ positioned at   the squares $(a_1, b_1)$, $(a_2, b_2)$ and $(r_m, c_m)$ would form a  NE-chain of length $3$ in $F$, which contradicts the hypothesis. This completes the proof. \qed

\begin{figure}[h!]
\begin{center}
\begin{picture}(100,30)
\setlength{\unitlength}{0.6mm}
\put( 0,0){\framebox(5,5) } \put(5,0){\framebox(5,5)  }  \put(10,0){\framebox(5,5)  }
 \put(15,0){\framebox(5,5) }  \put(20,0){\framebox(5,5) }  \put(25,0){\framebox(5,5) {$\bullet$}}
  \put(30,0){\framebox(5,5)  }  \put(35,0){\framebox(5,5)  }  \put(40,0){\framebox(5,5)  }

  \put( 0,5){\framebox(5,5) } \put(5,5){\framebox(5,5)  }  \put(10,5){\framebox(5,5)  }
 \put(15,5){\framebox(5,5) }  \put(20,5){\framebox(5,5) }  \put(25,5){\framebox(5,5)  }
  \put(30,5){\framebox(5,5) {$\bullet$} }  \put(35,5){\framebox(5,5)  }

\put( 0,10){\framebox(5,5) } \put(5,10){\framebox(5,5)  }  \put(10,10){\framebox(5,5)  }
 \put(15,10){\framebox(5,5) }  \put(20,10){\framebox(5,5) }  \put(25,10){\framebox(5,5){$\bullet$}  }
  \put(30,10){\framebox(5,5)  }

  \put( 0,15){\framebox(5,5) } \put(5,15){\framebox(5,5)  }  \put(10,15){\framebox(5,5)  }
 \put(15,15){\framebox(5,5) }  \put(20,15){\framebox(5,5){$\bullet$} }  \put(25,15){\framebox(5,5)  }

 \put( 0,20){\framebox(5,5) } \put(5,20){\framebox(5,5)  }  \put(10,20){\framebox(5,5)  }
 \put(15,20){\framebox(5,5){$\bullet$}  }  \put(20,20){\framebox(5,5)}

 \put( 0,25){\framebox(5,5) } \put(5,25){\framebox(5,5)  }  \put(10,25){\framebox(5,5)  }
 \put(15,25){\framebox(5,5){$\bullet$}  }

 \put( 0,30){\framebox(5,5) } \put(5,30){\framebox(5,5)  }  \put(10,30){\framebox(5,5) {$\bullet$} }
  \put( 0,35){\framebox(5,5) } \put(5,35){\framebox(5,5) {$\bullet$}  }
  \put( 0,40){\framebox(5,5){$\bullet$}  }

  \put(40,20){$\longrightarrow$}
  \put(40,25){$i=4$}
  \put( 60,0){\framebox(5,5) } \put(65,0){\framebox(5,5)  }  \put(70,0){\framebox(5,5)  }
 \put(75,0){\framebox(5,5) }  \put(80,0){\framebox(5,5) }  \put(85,0){\framebox(5,5) {$\bullet$}}
  \put(90,0){\framebox(5,5)  }  \put(95,0){\framebox(5,5)  }  \put(100,0){\framebox(5,5)  }

  \put( 60,5){\framebox(5,5) } \put(65,5){\framebox(5,5)  }  \put(70,5){\framebox(5,5)  }
 \put(75,5){\framebox(5,5) }  \put(80,5){\framebox(5,5) }  \put(85,5){\framebox(5,5)  }
  \put(90,5){\framebox(5,5) {$\bullet$} }  \put(95,5){\framebox(5,5)  }

\put( 60,10){\framebox(5,5) } \put(65,10){\framebox(5,5)  }  \put(70,10){\framebox(5,5)  }
 \put(75,10){\framebox(5,5) }  \put(80,10){\framebox(5,5) }  \put(85,10){\framebox(5,5){$\bullet$}  }
  \put(90,10){\framebox(5,5)  }

  \put( 60,15){\framebox(5,5) } \put(65,15){\framebox(5,5)  }  \put(70,15){\framebox(5,5)  }
 \put(75,15){\framebox(5,5) }  \put(80,15){\framebox(5,5){$\bullet$} }  \put(85,15){\framebox(5,5)  }

 \put( 60,20){\framebox(5,5) } \put(65,20){\framebox(5,5)  }  \put(70,20){\framebox(5,5)  }
 \put(75,20){\framebox(5,5){$\bullet$}  }  \put(80,20){\framebox(5,5)}

 \put( 60,25){\framebox(5,5) } \put(65,25){\framebox(5,5)  }  \put(70,25){\framebox(5,5)  }
 \put(75,25){\framebox(5,5)  }

 \put( 60,30){\framebox(5,5) } \put(65,30){\framebox(5,5)  }  \put(70,30){\framebox(5,5) {$\bullet$} }
  \put( 60,35){\framebox(5,5) } \put(65,35){\framebox(5,5) {$\bullet$}  }
  \put( 60,40){\framebox(5,5){$\bullet$}  }

\put(100,20){$\longrightarrow$}
\put(100,25){$i=6$}
  \put( 120,0){\framebox(5,5) } \put(125,0){\framebox(5,5)  }  \put(130,0){\framebox(5,5)  }
 \put(135,0){\framebox(5,5) }  \put(140,0){\framebox(5,5) }  \put(145,0){\framebox(5,5) {$\bullet$}}
  \put(150,0){\framebox(5,5)  }  \put(155,0){\framebox(5,5)  }  \put(160,0){\framebox(5,5)  }

  \put( 120,5){\framebox(5,5) } \put(125,5){\framebox(5,5)  }  \put(130,5){\framebox(5,5)  }
 \put(135,5){\framebox(5,5) }  \put(140,5){\framebox(5,5) }  \put(145,5){\framebox(5,5)  }
  \put(150,5){\framebox(5,5) {$\bullet$} }  \put(155,5){\framebox(5,5)  }

\put( 120,10){\framebox(5,5) } \put(125,10){\framebox(5,5)  }  \put(130,10){\framebox(5,5)  }
 \put(135,10){\framebox(5,5) }  \put(140,10){\framebox(5,5){$\bullet$}  }  \put(145,10){\framebox(5,5) }
  \put(150,10){\framebox(5,5)  }

  \put( 120,15){\framebox(5,5) } \put(125,15){\framebox(5,5)  }  \put(130,15){\framebox(5,5)  }
 \put(135,15){\framebox(5,5) }  \put(140,15){\framebox(5,5)  }  \put(145,15){\framebox(5,5)  }

 \put( 120,20){\framebox(5,5) } \put(125,20){\framebox(5,5)  }  \put(130,20){\framebox(5,5)  }
 \put(135,20){\framebox(5,5){$\bullet$}  }  \put(140,20){\framebox(5,5)}

 \put( 120,25){\framebox(5,5) } \put(125,25){\framebox(5,5)  }  \put(130,25){\framebox(5,5)  }
 \put(135,25){\framebox(5,5)  }

 \put( 120,30){\framebox(5,5) } \put(125,30){\framebox(5,5)  }  \put(130,30){\framebox(5,5) {$\bullet$} }
  \put( 120,35){\framebox(5,5) } \put(125,35){\framebox(5,5) {$\bullet$}  }
  \put( 120,40){\framebox(5,5){$\bullet$}  }

  \end{picture}
\end{center}
\caption{ An example of two iterations of $\alpha$ to a $01$-filling in $\mathcal{M}_3(9)$.  }\label{filling1}
\end{figure}
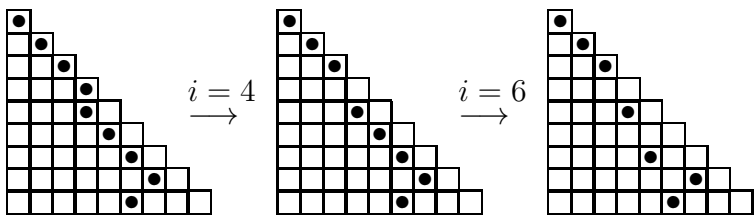

Lemma \ref{alpha0} states that the column $i$ that we find in the transformation $\alpha$ can only go right. Hence, there will be no  column containing at least two $1's$  in the resulting filling
after finitely many iterations of $\alpha$. Lemma \ref{alpha1} tells us that the resulting filling is a valid $01$-filling of $\Delta_n$ containing no NE-chain of length $3$. Therefore, we will get a $01$-filling in $N_3(n)$ after finitely applying  many iterations of $\alpha$ to a $01$-filling $F$ in $\mathcal{M}_3(n)$.
Define   $\phi(F)$ to be the resulting filling. Figure \ref{filling1} illustrates an example of   two iterations of $\alpha$ to a $01$-filling in $\mathcal{M}_3(9)$.

 {\noindent \bf The transformation $\beta$} \,\, Let $F$ be a valid filling of $\Delta_n$ which verifies  property (b1) and  contains no NE-chain of length $3$.  If every row contains a $1$ in $F$, then we simply define $\beta(F)=F$. Otherwise, find the lowest zero  row $i$. Suppose that the $1's$ below and  weakly to the left of  the square $(i,i)$ are positioned at the squares $(r_1, c_1)$, $(r_2, c_2), \ldots, (r_m, c_m)$ with $r_1<r_2<\ldots <r_m$.  Assume that $r_0=i$.
  Suppose that the topmost $1$ at column $i$ is positioned at the  square $(r_s,c_s)$.

     If there is at least one $1$ which is above and  to the right of the square $(r_s, c_s)$, then find the topmost square, say $(p,q)$, containing  such a $1$. Then we have $p=r_t+1$ for some $0\leq t\leq s-1$.
  Define $\beta(F)$ to be the $01$-filling of $\Delta_n$ obtained   from $F$ by the following procedure:
\begin{itemize}
 \item For all $1\leq  \ell \leq  t$,  replace the $1$ at the square    $(r_\ell, c_{\ell})$ with a $0$;
   \item  For all $0\leq \ell\leq t$,  fill  the square $(r_\ell, c_{\ell+1})$  with a $1 $ with the assumption $c_{t+1}=i$;

   \item Leave all the other squares fixed.
\end{itemize}

Otherwise,
   define $\beta(F)$ to be the $01$-filling of $\Delta_n$ obtained   from $F$ by the following procedure:
\begin{itemize}
 \item For all $1\leq \ell \leq m$, replace the $1$ at the square  $(r_\ell, c_\ell)$ with  a $0 $;
   \item  For all $0\leq \ell\leq m $,  fill  the square  $(r_\ell, c_{\ell+1})$  with a $1$ with the assumption $c_{m+1}=i$;
   \item Leave all the other squares fixed.
\end{itemize}

Now we  proceed to   show that the transformation  $\beta$   has the following   analogous properties of $\alpha$.
\begin{lemma}\label{beta0}
In $\beta(F)$, every row   below  row $i$ (including row $i$) contains  exactly  one $1$.
\end{lemma}
\pf It is obvious from the selection of row $i$ and the  construction of the transformation $\beta$. \qed

\begin{lemma}\label{beta1}
   The filling $\beta(F)$  is a valid  $01$-filling of $\Delta_n$  which verifies property (b1) and  contains no  NE-chain of length $3$.
 \end{lemma}
\pf It is obvious that $\beta(F)$  is a valid $01$-filling of $\Delta_n$ which verifies property (b1). Now we proceed to show that $\beta(F)$ contains no  NE-chain of length $3$. If not, suppose that the $1's$ positioned at the squares $(a_1, b_1)$, $(a_2, b_2)$ and $(a_3, b_3)$ form a NE-chain of length $3$, where $a_1<a_2<a_3$.

  Suppose that there is at least one $1$ which is above and to the right of the square $(r_s,c_s)$ in $F$. In this case,   all the squares  below row $r_t$ in $\beta(F)$ remain the same as those of $F$. Since  $F$ contains no NE-chain of length $3$, we must have $(a_1, b_1)=(r_t, i)$.   Hence, the $1's$ positioned at the squares $(p,q)$, $(a_2, b_2)$ and $(a_3, b_3)$ form a NE-chain of length $3$, which contradicts the hypothesis.
Thus, $F$ does not contain a $1$ which is above and to the right of the square $(r_s,c_s)$. According to the construction of $\beta(F)$, one of $(a_1, b_1)$, $(a_2, b_2)$   and $(a_3, b_3)$ must fall in $(r_m,i)$. Since there is no $1$ which is below row $r_m$ and to the left  of column $i$ in $\beta(F)$, we have $(a_3, b_3)=(r_m,i)$. Then the $1's$ positioned at the squares $(a_1, b_1)$, $(a_2, b_2)$ and $(r_m, i)$ would form a   NE-chain of length $3$ in $F$, which contradicts the hypothesis. This completes the proof. \qed

Lemma \ref{beta0} states that the row $i$ that we find in the transformation $\beta$ can only go upside. Hence,  there will be no zero row
after finitely many iterations of $\beta$. Lemma \ref{beta1} tells us that the resulting filling is a valid $01$-filling of $\Delta_n$ containing no NE-chain of length $3$. Hence, we will get a $01$-filling in $\mathcal{M}_3(n)$ after finitely applying  many iterations of $\beta$ to a $01$-filling $F$ in $\mathcal{N}_3(n)$.
Define   $\psi(F)$ to be the resulting filling.

\begin{theorem}\label{bijection}
The maps $\phi$ and $\psi$ induce a bijection between the set $\mathcal{M}_3(n)$ and the set $\mathcal{N}_3(n)$.
\end{theorem}
\pf It suffices to show that the maps $\phi$ and $\psi$ are inverses of each other. First, we proceed to show that $\phi$ is the inverse of the map $\psi$, that is, $\phi(\psi(F))=F$ for any   $01$-filling $F\in \mathcal{N}_3(n)$. To this end, it suffices to show that $\alpha(\beta^k(F))=\beta^{k-1}(F)$.     Suppose that at the $k$th application of $\beta$ to $\beta^{k-1}(F)$, the selected row is row $i$. Suppose that the $1's$ below and  weakly to the left of  the square $(i,i)$ are positioned at the squares $(r_1, c_1)$, $(r_2, c_2), \ldots, (r_m, c_m)$ with $r_1<r_2<\ldots <r_m$.  Assume that $r_0=i$.
  Suppose that the topmost  $1$ at column $i$ is positioned at the  square $(r_s,c_s)$. We have two cases.

If   there is at least one 1 which is above and to the right of the square $(r_s,c_s)$ in $\beta^{k-1}(F)$, then
find the topmost square $(p,q)$ containing such a 1. Assume that $p=r_t+1$ for some $0\leq t<s$.
  From the construction of the transformation $\beta$,  the square  $(r_\ell, c_{\ell+1})$ of $\beta^{k}(F)$  is filled with a $1$ for all $1\leq \ell\leq t $  with the assumption $c_{t+1}=i$ and all the other squares remain the same as those of $\beta^{k-1}(F)$. Clearly, in $\beta^{k}(F)$,  all the columns to the left of column  $i$ contains at most one $1$, and  column $i$ contains exactly two $1's$. Hence, when we apply the transformation $\alpha$ to $\beta^{k}(F)$, the column that we select is just    column     $i$ and the topmost $1$   at column $i$ is positioned at the square $(r_t, i)$. Moreover,  the $1$ positioned at the square $(p,q)$ is to the right of square $(r_t, i)$ and $p=r_t+1$.  From the construction of $\alpha$, it is not dificult to check that $\alpha(\beta^{k}(F))=\beta^{k-1}(F)$.

If   there  does not exists any 1 which is above and to the right of the square $(r_s, c_s)$ in $\beta^{k-1}(F)$,
  From the construction of the transformation $\beta$,  the square  $(r_\ell, c_{\ell+1})$ of $\beta^{k}(F)$  is filled with a $1$ for all $1\leq \ell\leq m $  with the assumption $r_{m+1}=i$  and all the other squares remain the same as those of $\beta^{k-1}(F)$. Clearly, in $\beta^{k}(F)$,  all the columns to the left of column  $i$ contains at most one $1$,  in which column $i$ contains exactly two $1's$.  Hence, when we apply the transformation $\alpha$ to $\beta^{k}(F)$, the column that we select is just    column     $i$ and the topmost $1$   at column $i$ is positioned at the square $(r_{s-1}, c_s)$.
Notice that  there does not exist any $1$ which is above and to the right of square $(r_s,c_s)$ in $\beta^{k-1}(F)$.
Hence,  there is no $1$ at row $r_{s-1}+1$ which  is  to the right  of the square $(r_{s-1}, c_s)$.  From the construction of $\alpha$, it is easily seen that $\alpha(\beta^{k}(F))=\beta^{k-1}(F)$.

Combining the two above cases,  we have deduced that   $\alpha(\beta^k(F))=\beta^{k-1}(F)$. By similar arguments, one can verify that $\beta(\alpha^k(F))=\alpha^{k-1}(F)$ for any $01$-filling $F$ of $\mathcal{M}_3(n)$. The details are omitted here.   Hence, the maps $\phi$ and 	$\psi$ are inverses of each other. Thus,
the maps $\phi$ and $\psi$  induce a bijection  between the set $\mathcal{M}_3(n)$ and the set $\mathcal{N}_3(n)$ as claimed.
 \qed

Combining Theorems \ref{th1}, \ref{th2} and \ref{bijection}, we are led to a bijective proof of Conjecture \ref{con}.

\section{Bijective proof of Conjecture \ref{Yan} }
In this section, we shall give a new bijective proof of Conjecture \ref{Yan} relying on the bijection $\phi$.

In the following,  a $01$-filling in $\mathcal{M}_3(n)$   will be identified with a sequence $\{( 1, a_1),$ $  (2, a_2), $ $ \cdots, (n, a_n)\}$, where $1\leq a_i\leq i$ and $a_i=k$ if and only if  there is a $1$ in the $i$th row and $k$th column.
In the course of proving Conjecture \ref{Yan}, Yan \cite{Yan} provided a bijection $\gamma$ between the set $\mathcal{PA}_3(n+1)$ and the set $\mathcal{M}_3(n)$. Let $x=x_1x_2\cdots x_{n+1}\in \mathcal{PA}_{3}(n+1)$.  Define $\gamma(x)=\{( 1, a_1),  (2, a_2), $ $\cdots$  ,$(n, a_n)\}$  where $a_i=i+x_{i+1}-\as(x_1x_2\cdots x_{i+1})$ for all $i=1, 2, \cdots, n$. For example,  let $x=012340415\in \mathcal{PA}_{3}(9)$. Then we have
 $$\gamma(x)=\{(1, 1), (2, 2), (3,3), (4,4), (5,1), (6, 5), (7, 3), (8, 7) \}. $$

The inverse of the map $\gamma$ is defined as follows. Let $F=\{( 1, a_1),  (2, a_2),  \cdots  ,(n, a_n)\}$. Define $\gamma^{-1}(F)=(x_1,x_2,\cdots, x_{n+1})$ inductively  as follows:
  \begin{itemize}
  \item   $x_1=0$ and $x_2=1$;
   \item  if $a_{i-1 }<  a_{i }$, then $x_{i+1}=\as(x_1x_2\cdots x_{i})+1+a_{i}-i$  for all $2\leq i\leq n$ ;
 \item if  $a_{i-1}\geq  a_{i}$, then $x_{i+1}=\as(x_1x_2\cdots x_{i})+a_{i}-i$  for all $2\leq i\leq n$.

\end{itemize}

Recall that  Krattenthaler \cite{kra} also established      a bijection  between set partitions of $[n+1]$ and $01$-fillings of $\Delta_{n}$  in which
   every row and every column contain  at most one $1$.
 Given a set partition $\pi$ of $[n]$, we can get a $01$-filling of $\Delta_{n}$ by putting a $1$ in the square $(j-1,i)$  if  $(i,j)$ is an arc in its  linear representation.   From the construction of Krattenthaler's bijection,    a   $k$-nesting of a set partition  corresponds to a NE-chain of length $k$ in its corresponding $01$-filling.

 Denote by $\mathcal{P}_k(n)$ the set of $01$-fillings of $\Delta_{n}$  in which every row and every column contain  at most one $1$, and there is no NE-chain of length $k$.

 The following result follows immediately from  Krattenthaler's bijection \cite{kra}.
 \begin{theorem}\label{th4}
   There is a one-to-one correspondence between    the set $\mathcal{C}_k(n+1)$ and the set $\mathcal{P}_k(n)$.
\end{theorem}
 By Theorem \ref{th4}, in order to provide a bijection between $\mathcal{A}_3(n)$ and $C_3(n)$, it suffices to establish a bijection between the set $\mathcal{A}_3(n+1)$ and the set $\mathcal{P}_3(n)$.

 In a $01$-filling, if both row $i$ and column $i$ are zero,  then row (column) $i$  is said to be {\em critical}.

\begin{theorem}\label{mainth}
There is a bijection between the set $\mathcal{A}_3(n+1)$ and the set $\mathcal{P}_3(n)$.
\end{theorem}
\pf  First we shall describe a map $\delta$ from  the set $\mathcal{A}_3(n+1)$ to the set $\mathcal{P}_3(n)$.  Let $x\in \mathcal{A}_3(n+1) $.   It is apparent that the ascent  sequence
 $x$ can be written   as $x_1^{c_1} x_2^{c_2} \cdots  x_{k+1}^{c_{k+1}}$, where $x_i\neq x_{i+1}$  and $c_i\geq 1$ for  all $i\geq 1$. Let $x'=x_1 x_2\cdots  x_{k+1}$. Obviously, $x'$ is a  primitive ascent sequence in $\mathcal{PA}_3(k+1)$. Let $F=\gamma(x')$ and $F'=\phi(F)$. Clearly, we have $F\in \mathcal{M}_3(k)$ and $F'\in \mathcal{N}_3(k)$.  Now we can generate a $01$-filling $F''$ of $\Delta_{n}$ from $F'$ by inserting $c_1-1$ consecutive  zero rows immediately above row $1$ and    $c_1-1$ consecutive  zero  columns immediately to the left of column $1$,  and inserting $c_i$  consecutive  zero rows immediately below  row $i$  and       $c_i$  consecutive zero columns immediately to the right of column $i$ for all $1\leq i\leq k$. Define $\delta(x)=F''$.  It is not difficult to see that  the resulting filling $F''$ is an element of $\mathcal{P}_3(n)$.  This implies that the map $\delta$ is well defined.

 In order to prove that $\delta$ is a bijection, we construct a map $\delta'$ from the set $\mathcal{P}_3(n)$ to the set $\mathcal{A}_3(n+1)$.
 Given a   $01$-filling  $F\in \mathcal{P}_3(n)$, we can recover an ascent sequence $\delta'(F)$ as follows. Suppose that  there are $k$  non-critical rows in $F$.  Let rows $i_1$, $i_2$, $\ldots$, $i_k$  be  the non-critical rows of $F$.   Assume that there are $c_1$ critical rows immediately above row $i_1$, and $c_{\ell+1}$ critical rows immediately below row $i_\ell$ for all $1\leq \ell\leq k$. Denote  $F'$ the $01$-filling obtained from $F$ by removing all the  critical rows and columns from $F$. Moreover, let $F''=\psi(F')$ and $x=x_1x_2\ldots x_{k+1}=\gamma^{-1}(F'')$. It is easily seen that $F'\in \mathcal{N}_3(k)$, $F''\in \mathcal{M}_3(k)$ and $x\in \mathcal{PA}_3(k+1)$. Let $\delta'(F)=x_1^{c_1+1}x_2^{c_{2}}\ldots x_{k+1}^{c_{k+1}} $. It is apparent  that we have $\delta'(F)\in \mathcal{A}_{3}(n+1)$.

  Property  $(b1)$ ensures that the
inserted rows and columns  in the construction of $\delta$ are exactly the removed rows and columns
  in the construction of $\delta'$.  Thus the map $\delta'$ is the inverse of
the map $\delta$. This implies that $\delta$ is bijection.  \qed

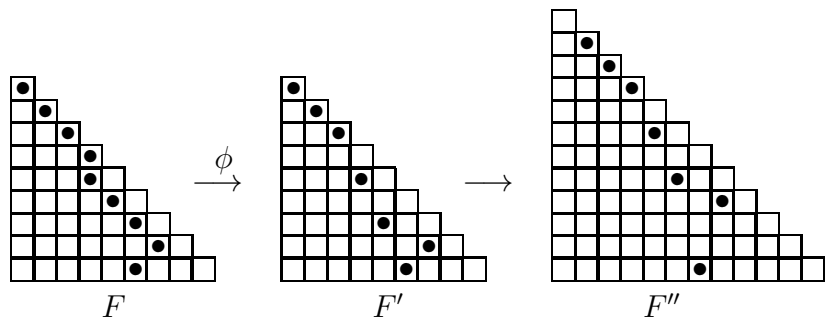
\begin{figure}[h!]
\begin{center}
\begin{picture}(100,45)
\setlength{\unitlength}{0.6mm}
\put( 0,0){\framebox(5,5) } \put(5,0){\framebox(5,5)  }  \put(10,0){\framebox(5,5)  }
 \put(15,0){\framebox(5,5) }  \put(20,0){\framebox(5,5) }  \put(25,0){\framebox(5,5) {$\bullet$}}
  \put(30,0){\framebox(5,5)  }  \put(35,0){\framebox(5,5)  }  \put(40,0){\framebox(5,5)  }

  \put( 0,5){\framebox(5,5) } \put(5,5){\framebox(5,5)  }  \put(10,5){\framebox(5,5)  }
 \put(15,5){\framebox(5,5) }  \put(20,5){\framebox(5,5) }  \put(25,5){\framebox(5,5)  }
  \put(30,5){\framebox(5,5) {$\bullet$} }  \put(35,5){\framebox(5,5)  }

\put( 0,10){\framebox(5,5) } \put(5,10){\framebox(5,5)  }  \put(10,10){\framebox(5,5)  }
 \put(15,10){\framebox(5,5) }  \put(20,10){\framebox(5,5) }  \put(25,10){\framebox(5,5){$\bullet$}  }
  \put(30,10){\framebox(5,5)  }

  \put( 0,15){\framebox(5,5) } \put(5,15){\framebox(5,5)  }  \put(10,15){\framebox(5,5)  }
 \put(15,15){\framebox(5,5) }  \put(20,15){\framebox(5,5){$\bullet$} }  \put(25,15){\framebox(5,5)  }

 \put( 0,20){\framebox(5,5) } \put(5,20){\framebox(5,5)  }  \put(10,20){\framebox(5,5)  }
 \put(15,20){\framebox(5,5){$\bullet$}  }  \put(20,20){\framebox(5,5)}

 \put( 0,25){\framebox(5,5) } \put(5,25){\framebox(5,5)  }  \put(10,25){\framebox(5,5)  }
 \put(15,25){\framebox(5,5){$\bullet$}  }

 \put( 0,30){\framebox(5,5) } \put(5,30){\framebox(5,5)  }  \put(10,30){\framebox(5,5) {$\bullet$} }
  \put( 0,35){\framebox(5,5) } \put(5,35){\framebox(5,5) {$\bullet$}  }
  \put( 0,40){\framebox(5,5){$\bullet$}  }
\put(20,-8){$F$}
\put(40,20){$\longrightarrow$}
\put(45,25){$\phi$}
  \put( 60,0){\framebox(5,5) } \put(65,0){\framebox(5,5)  }  \put(70,0){\framebox(5,5)  }
 \put(75,0){\framebox(5,5) }  \put(80,0){\framebox(5,5) }  \put(85,0){\framebox(5,5) {$\bullet$}}
  \put(90,0){\framebox(5,5)  }  \put(95,0){\framebox(5,5)  }  \put(100,0){\framebox(5,5)  }

  \put( 60,5){\framebox(5,5) } \put(65,5){\framebox(5,5)  }  \put(70,5){\framebox(5,5)  }
 \put(75,5){\framebox(5,5) }  \put(80,5){\framebox(5,5) }  \put(85,5){\framebox(5,5)  }
  \put(90,5){\framebox(5,5) {$\bullet$} }  \put(95,5){\framebox(5,5)  }

\put( 60,10){\framebox(5,5) } \put(65,10){\framebox(5,5)  }  \put(70,10){\framebox(5,5)  }
 \put(75,10){\framebox(5,5) }  \put(80,10){\framebox(5,5){$\bullet$}  }  \put(85,10){\framebox(5,5) }
  \put(90,10){\framebox(5,5)  }

  \put( 60,15){\framebox(5,5) } \put(65,15){\framebox(5,5)  }  \put(70,15){\framebox(5,5)  }
 \put(75,15){\framebox(5,5) }  \put(80,15){\framebox(5,5)  }  \put(85,15){\framebox(5,5)  }

 \put( 60,20){\framebox(5,5) } \put(65,20){\framebox(5,5)  }  \put(70,20){\framebox(5,5)  }
 \put(75,20){\framebox(5,5){$\bullet$}  }  \put(80,20){\framebox(5,5)}

 \put( 60,25){\framebox(5,5) } \put(65,25){\framebox(5,5)  }  \put(70,25){\framebox(5,5)  }
 \put(75,25){\framebox(5,5)  }

 \put( 60,30){\framebox(5,5) } \put(65,30){\framebox(5,5)  }  \put(70,30){\framebox(5,5) {$\bullet$} }
  \put( 60,35){\framebox(5,5) } \put(65,35){\framebox(5,5) {$\bullet$}  }
  \put( 60,40){\framebox(5,5){$\bullet$}  }
 \put(80,-8){$F'$}

\put(100,20){$\longrightarrow$}

  \put( 120,0){\framebox(5,5) }  \put( 125,0){\framebox(5,5) }  \put( 130,0){\framebox(5,5) }  \put( 135,0){\framebox(5,5) }
 \put( 140,0){\framebox(5,5) }  \put( 145,0){\framebox(5,5) }  \put( 150,0){\framebox(5,5){$\bullet$} }  \put( 155,0){\framebox(5,5) }
  \put( 160,0){\framebox(5,5) }  \put( 165,0){\framebox(5,5) }  \put( 170,0){\framebox(5,5) }  \put( 175,0){\framebox(5,5) }

  \put( 120,5){\framebox(5,5) }  \put( 125,5){\framebox(5,5) }  \put( 130,5){\framebox(5,5) }  \put( 135,5){\framebox(5,5) }
 \put( 140,5){\framebox(5,5) }  \put( 145,5){\framebox(5,5) }  \put( 150,5){\framebox(5,5)  }  \put( 155,5){\framebox(5,5) }
  \put( 160,5){\framebox(5,5) }  \put( 165,5){\framebox(5,5) }  \put( 170,5){\framebox(5,5) }

  \put( 120,10){\framebox(5,5) }  \put( 125,10){\framebox(5,5) }  \put( 130,10){\framebox(5,5) }  \put( 135,10){\framebox(5,5) }
 \put( 140,10){\framebox(5,5) }  \put( 145,10){\framebox(5,5) }  \put( 150,10){\framebox(5,5)  }  \put( 155,10){\framebox(5,5) }
  \put( 160,10){\framebox(5,5) }  \put( 165,10){\framebox(5,5) }

  \put( 120,15){\framebox(5,5) }  \put( 125,15){\framebox(5,5) }  \put( 130,15){\framebox(5,5) }  \put( 135,15){\framebox(5,5) }
 \put( 140,15){\framebox(5,5) }  \put( 145,15){\framebox(5,5) }  \put( 150,15){\framebox(5,5)  }  \put( 155,15){\framebox(5,5) {$\bullet$}}
  \put( 160,15){\framebox(5,5) }

  \put( 120,20){\framebox(5,5) }  \put( 125,20){\framebox(5,5) }  \put( 130,20){\framebox(5,5) }  \put( 135,20){\framebox(5,5) }
 \put( 140,20){\framebox(5,5) }  \put( 145,20){\framebox(5,5){$\bullet$}  }  \put( 150,20){\framebox(5,5)  }  \put( 155,20){\framebox(5,5) }

 \put( 120,25){\framebox(5,5) }  \put( 125,25){\framebox(5,5) }  \put( 130,25){\framebox(5,5) }  \put( 135,25){\framebox(5,5) }
 \put( 140,25){\framebox(5,5) }  \put( 145,25){\framebox(5,5)   }  \put( 150,25){\framebox(5,5)  }

  \put( 120,30){\framebox(5,5) }  \put( 125,30){\framebox(5,5) }  \put( 130,30){\framebox(5,5) }  \put( 135,30){\framebox(5,5) }
 \put( 140,30){\framebox(5,5){$\bullet$} }  \put( 145,30){\framebox(5,5)   }

 \put( 120,35){\framebox(5,5) }  \put( 125,35){\framebox(5,5) }  \put( 130,35){\framebox(5,5) }  \put( 135,35){\framebox(5,5) }
 \put( 140,35){\framebox(5,5)  }

  \put( 120,40){\framebox(5,5) }  \put( 125,40){\framebox(5,5) }  \put( 130,40){\framebox(5,5) }  \put( 135,40){\framebox(5,5) {$\bullet$}}

  \put( 120,45){\framebox(5,5) }  \put( 125,45){\framebox(5,5) }  \put( 130,45){\framebox(5,5){$\bullet$} }
   \put( 120,50){\framebox(5,5) }  \put( 125,50){\framebox(5,5)     {$\bullet$} }
   \put( 120,55){\framebox(5,5) }
   \put(140,-8){$F''$}

  \end{picture}
\end{center}
\caption{ A $01$-filling $F\in \mathcal{M}_3(9)$, a $01$-filling $F'\in \mathcal{N}_3(9)$ and a $01$-filling $F''\in \mathcal{P}_3(12)$.}\label{filling2}
\end{figure}

For example, let $x=001234345664$ be   an   ascent sequence  in $\mathcal{A}_3(13)$. Then $x$ can be written as $0^21^12^13^14^13^14^15^16^{2}4  $. Let $x'=0123434564$, which is an element of $\mathcal{PA}_3(10)$.
By applying the map $\gamma$ to $x'$, we get a $01$-filling
$$F=\gamma(x')=\{(1,1), (2,2,), (3,3), (4,4), (5,4), (6,5), (7,6), (8,7), (9,6)\}\in \mathcal{M}_3(9) $$  illustrated  in Figure \ref{filling2}.  Then by applying the map $\phi$ to $F$, we get a $01$-filling $F'$ as shown in Figure \ref{filling2}. Finally, we obtain a $01$-filling $F''\in \mathcal{P}_3(12)$ by adding   one zero row immediately above row $1$, one zero column immediately  to the left of column $1$, two consecutive zero rows immediately below row $8$,  and two consecutive zero columns immediately  to the right of column $8$, see Figure \ref{filling2}.

  Combining  Theorems  \ref{th4}  and \ref{mainth}, we get    a new bijective  proof of  Conjecture \ref{Yan}.

\noindent{\bf Acknowledgments.}
 This work was supported by  the National Natural Science Foundation of China (11571320 and 11671366) and  Zhejiang Provincial Natural Science Foundation of China ( LY15A010008).


\end{document}